\documentclass[english, a4paper]{article}
\usepackage[T1]{fontenc}
\usepackage[cp1251]{inputenc}
\usepackage{amssymb}
\usepackage{babel}

\usepackage{amsfonts}

\usepackage{amsmath, amsthm, amscd, graphicx,amssymb}
\newtheorem{rem}{Remark}
\newtheorem{teo}{Theorem}
\newcommand{\dd}{\operatorname{d}\nolimits}

\def\ksi{\xi}
\def\h{u}
\def\AA{\mathcal{A}}
\def\G{\mathcal{G}}
\def\R{{\mathbb R}}
\newcommand{\Id}{e}
\newcommand{\SO}{\operatorname{SO(3)}\nolimits}
\newcommand{\SE}{\operatorname{SE(3)}\nolimits}
\newcommand{\eexp}{\operatorname{exp}\nolimits}
\newcommand{\ee}{\operatorname{e}\nolimits}
\newcommand{\cn}{\operatorname{cn}\nolimits}
\newcommand{\sn}{\operatorname{sn}\nolimits}
\newcommand{\dn}{\operatorname{dn}\nolimits}
\newcommand{\sh}{\operatorname{sinh}\nolimits}

\newcommand{\tal}{\alpha}
\newcommand{\tbe}{\beta}
\newcommand{\tga}{\theta}
\begin{document}

\title{Extremal Controls in the Sub-Riemannian Problem on the Group of Motions of Euclidean Space\footnote{This work is supported by the Russian Science Foundation
under grant 17-11-01387 and performed in Ailamazyan Program Systems
Institute of Russian Academy of Sciences.}
}

\author{A.~Mashtakov, A.~Popov \footnote{CPRC, Ailamazyan Program Systems Institute of RAS}}

\maketitle

\begin{abstract}
For the sub-Riemannian problem on the group of motions of Euclidean space we present explicit formulas for extremal controls in a special case, when one of the initial momenta is fixed.
\end{abstract}
%%%%%%%%%%%%%%%%%%%%%%%%%%%%%%%%%%%%%%%%%%%%%%%%%%%%%%%%%%%%%%%%%%%%%%%%%%%%%%%%%%%%%%%%%%%%%%%%%%%%%%%%%%%%%%%%%%%%%
\begin{flushright}
\textit{To the memory of Vladimir Arnold}
\end{flushright}
%%%%%%%%%%%%%%%%%%%%%%%%%%%%%%%%%%%%%%%%%%%%%%%%%%%%%%%%%%%%%%%%%%%%%%%%%%%%%%%%%%%%%%%%%%%%%%%%%%%%%%%%%%%%%%%%%%%%%
%%%%%%%%%%%%%%%%%%%%%%%%%%%%%%%%%%%%%%%%%%%%%%%%%%%%%%%%%%%%%%%%%%%%%%%%%%
\section{Introduction}\label{sec:1}
In this paper, we consider a sub-Riemannian (SR) problem on the group of motions of Euclidean space $\SE$. It can be interpreted as a problem of optimal motion of a rigid body in $\R^3$ with nonintegrable constraints~\cite{6p5}.  Solution curves to the problem have applications in image processing (tracking of neural fibres and blood vessels in MRI and CT images of human brain); and in robotics (motion planing problem for an aircraft, that can move forward/backward).

The sub-Riemannian problem on $\SE$ can be seen as follows. By given two orthonormal frames $N_0 = \{v_0^1, v_0^2, v_0^3\}$ and $N_1 = \{v_1^1, v_1^2, v_1^3\}$ attached respectively at two given points $q_0 = (x_0, y_0, z_0)$ and $q_1 = (x_1, y_1, z_1)$ in space $\R^3$, to find an optimal motion that transfers $q_0$ to $q_1$ such that the frame $N_0$ is transferred to the frame $N_1$. The frame can move forward or backward along one of the vector chosen in the frame and rotate simultaneously via two (of three) prescribed axes. The required motion should be optimal in the sense of minimal length in the space $\SE \cong \R^3 \times \SO$.  

The two-dimensional analog of this problem was studied as a possible model of the mechanism used by the visual cortex V1 of the human brain to reconstruct curves that are partially corrupted or hidden from observation. The two-dimensional model was initially due to~\cite{Citti_Sarti_Petitot}, where the authors recognized the sub-Riemannian Euclidean motion group structure of the problem. The related SR problem in $\operatorname{SE(2)}$ was solved in~\cite{SachkovSE2}, where in particular explicit formulas for the geodesics have been derived in SR arclength parameterization. Later, an alternative expression in spatial arclength parameterization for cuspless SR geodesics was derived in~\cite{DuitsSE2}. Application to contour completion in corrupted images was studied in~\cite{NMTMA}. The problem was also studied in~\cite{Hladky2010}.
However, many imaging applications such as diffusion weighted magnetic resonance imaging (DW-MRI) require an extension to three dimensions~\cite{Duits1,Duits2}, which motivates us to study the problem on $\SE$.

The Lie group $\SE$ of Euclidean motions of space $\R^3$ is generated by translations and rotations about coordinate axes. It is 
parameterized by matrices
\begin{eqnarray}
\label{eq:param}
\left(\begin{array}{cccc}
 \cos \tal  \cos \tbe  & -\cos \tbe  \sin \tal  & \sin \tbe  & x \\
 \cos \tga  \sin \tal +\cos \tal  \sin \tbe  \sin \tga  & \cos \tal  \cos \tga -\sin \tal  \sin \tbe  \sin \tga  & -\cos \tbe \sin \tga  & y \\
 \sin \tal  \sin \tga -\cos \tal  \cos \tga  \sin \tbe  & \cos \tga  \sin \tal  \sin \tbe +\cos \tal  \sin \tga  & \cos \tbe \cos \tga  & z \\
 0 & 0 & 0 & 1 \\
\end{array}\right),
\end{eqnarray}
where $\tga \in [-\frac{\pi }{2},\frac{\pi }{2})$, $\tbe \in[-\pi ,\pi )$, $\tal \in[0,2 \pi )$ are angles of rotation about the axes $OX$, $OY$, $OZ$; 
 and $(x,y,z) \in \R^3$ are coordinates with respect to the axes. 

Let us choose 
\begin{eqnarray*}
&& \AA_1=\cos{\tal}\cos{\tbe}\,\partial_x + (\sin{\tal}\cos{\tga} + \cos{\tal}\sin{\tbe}\sin{\tga})\,\partial_y + (\sin{\tal}\sin{\tga} -\cos{\tal}\sin{\tbe}\cos{\tga})\,\partial_z,\\
&& \AA_2=-\sin{\tal}\cos{\tbe}\,\partial_x + (\cos{\tal}\cos{\tga} - \sin{\tal}\sin{\tbe}\sin{\tga})\,\partial_y + (\cos{\tal}\sin{\tga} + \sin{\tal}\sin{\tbe}\cos{\tga})\,\partial_z,\\
&& \AA_3=\sin{\tbe}\,\partial_x - \cos{\tbe}\sin{\tga}\,\partial_y +  \cos{\tbe}\cos{\tga}\,\partial_z,\\
&& \AA_4= -\cos{\tal}\tan{\tbe}\,\partial_{\tal} + \sin{\tal}\,\partial_{\tbe} +  \cos{\tal}\sec{\tbe}\,\partial_{\tga},\\
&& \AA_5= \sin{\tal}\tan{\tbe}\,\partial_{\tal} + \cos{\tal}\,\partial_{\tbe} -  \sin{\tal}\sec{\tbe}\,\partial_{\tga},\\
&& \AA_6= \partial_{\tal}
\end{eqnarray*}
as the basis left-invariant vector fields agreed with parameterization~(\ref{eq:param}).

We consider the sub-Riemannian (SR) manifold $(\SE,\Delta,\G_{\ksi})$, see~\cite{1}. Here
$\Delta$ is a left-invariant distribution generated by the vector fields $\AA_3,\AA_4,\AA_5$;
$\G_{\ksi}$ is an inner product on $\Delta$ defined by
\begin{eqnarray*}
&& \G_{\ksi} = \ksi^2 \omega^3 \otimes \omega^3 +\omega^4 \otimes \omega^4 + \omega^5 \otimes \omega^5,
\end{eqnarray*}
with  $\ksi>0$ is a constant and $\omega^i$ are basis one forms satisfying
$$\langle \omega^i, \AA^j \rangle = \delta_j^i, \qquad \delta_i^j = 0, \text{ if } i\neq j, \quad \delta_i^i = 1.$$

We study a problem of finding sub-Riemannian length minimizers: By given boundary conditions, to find a Lipschitzian curve $\gamma : [0,t_1]\rightarrow \SE$, such that $\dot{\gamma}(t) \in \Delta$ for almost all $t \in (0, t_1)$ and $\gamma$ minimizes a functional of sub-Riemannian length
$$l(\gamma) = \int_0^{t_1} \sqrt{\G_\xi(\dot{\gamma}(t), \dot{\gamma}(t))} \dd t.$$

SR geodesics  are curves in $\SE$ whose sufficiently short arcs are SR minimizers. They satisfy the Pontryagin maximum principle, and the corresponding controls are called extremal controls.

Due to left-invariance of the problem one can fix the initial value $\gamma(0) = \Id$, where $\Id$ is the identical transformation of $\R^3$. Then the sub-Riemannian problem is equivalent to the following optimal control problem~\cite{1,2}:
\begin{eqnarray*}
&\dot \gamma = u_3 \AA_3 + u_4 \AA_4 + u_5 \AA_5, \label{eq:se3sys} \\
&\gamma(0) = \Id, \qquad \gamma(t_1) = q, \label{eq:se3bounds}\\
&l(\gamma) = \int_{0}^{t_1} \sqrt{\ksi^2 u_3^2\left(t\right) +  u_4^2\left(t\right) + u_5^2\left(t\right)}\, \dd t \to \min, \label{intsrlength}\\
\end{eqnarray*}
where the controls $u_3$, $u_4$, $u_5$ are real valued functions from $L_{\infty}(0,t_1)$.

The Cauchy-Schwarz inequality implies that the minimization problem for the sub-Riemannian length functional $l$ is equivalent to the minimization problem for the action functional
\begin{equation*}\label{eq:actionint}
J(\gamma) = \frac12 \int_0^{t_1} \left(\ksi^2 u_3^2\left(t\right) +u_4^2\left(t\right) + u_5^2\left(t\right)\right) \dd t \to \min,
\end{equation*}
with fixed $t_1>0$.

In paper~\cite{7}, the authors show that the problem can be reduced to the case $\xi=1$, and that application of the Pontryagin maximum principle leads to the following Hamiltonian system:
\begin{equation} \label{eq:hamsys}
\begin{array}{ll}
\begin{cases}
\dot{\h}_1 = -\h_3 \h_5,\\
\dot{\h}_2 = \h_3 \h_4,\\
\dot{\h}_3 = \h_1 \h_5 - \h_2 \h_4,\\
\dot{\h}_4 = \h_2 \h_3 - \h_5 \h_6,\\
\dot{\h}_5 = \h_4 \h_6 - \h_1 \h_3,\\
\dot{\h}_6 = 0,
\end{cases} %\label{vertpart}
& \begin{cases}
\dot{x} = \h_3 \sin\tbe, \\
\dot{y} = -\h_3\cos\tbe \sin\tga, \\
\dot{z} = \h_3\cos\tbe \cos\tga,\\
\dot{\tga} = \sec\tbe (\h_4  \cos\tal - \h_5 \sin\tal), \\
\dot{\tbe} = \h_4 \sin\tal +\h_5\cos\tal, \\
\dot{\tal} = -(\h_4 \cos\tal - \h_5 \sin\tal) \tan\tbe,
\end{cases} \\ %\label{horpart} \\
 \text{--- the vertical part (for extremal controls),}
&  \text{--- the horizontal part (for geodesics).}
\end{array}
\end{equation}

The vertical part describes dynamics of the extremal controls $u_3$, $u_4$, $u_5$ together with the remaining momentum components $u_1$, $u_2$, $u_6$. SR geodesics are solutions to the horizontal part. 

In this paper we focus on the simplest case $u_6 = 0$, as the most important for applications, in particular, for tracking of neural fibres and blood vessels in MRI and CT images of human brain~\cite{7}. In this case, the system on extremal controls becomes
\begin{equation} \label{eq:extcontsys}
\dot{\h}_1 = -\h_3 \h_5,\;\;
\dot{\h}_2 = \h_3 \h_4,\;\;
\dot{\h}_3 = \h_1 \h_5 - \h_2 \h_4,\;\;
\dot{\h}_4 = \h_2 \h_3,\;\;
\dot{\h}_5 = -\h_1 \h_3.\;\;
\end{equation}

We generalize results of~\cite{7}, where, in particular, the extremal controls are found in the case when the geodesics do not have cusps in their spatial projection. Such geodesics admit parametrization by spatial arclength, which leads to expression for the extremal controls in elementary functions.  Now, we relax the 'cuspless' assumption and derive explicit expression for $u_1, \ldots, u_5$ in terms of Jacobi elliptic functions.    

In Section~\ref{sec:2} we show, that if the function $u_3$ is known, then the first, the second, the fourth and the fifth equations of system~(\ref{eq:extcontsys}) allow us to express $u_k$, $k \in \{1,2,4,5\}$ via the initial values $u_k(0)$. Then by substitution of $u_k$ in the third equation of system~(\ref{eq:extcontsys}) we obtain an ordinary differential equation on $u_3$. Solution to this equation is presented in Section~\ref{sec:3}. 

\begin{rem}
Finding a parameterization of SR geodesics is a nontrivial problem. First natural question arises on a theoretical possibility of such parameterization in some reasonable sense --- the question of integrability of the Hamiltonian system. It was shown in~\cite[Thm.~2]{7}, that~(\ref{eq:hamsys}) is Liouville integrable, since it has a complete set of functionally independent first integrals in involution: $\h_6$, the Hamiltonian $H=\frac{1}{2}(\h_3^2+\h_{4}^2+\h_{5}^2)$, a Casimir function $W=\h_{1}\h_4+\h_{2}\h_{5}+\h_3\h_6$, and the right-invariant Hamiltonians  %$\rho_1$, $\rho_2$, $\rho_3$, given by 
\begin{equation*}
\begin{array}{l}
\rho_{1}  = -\h_1 \cos\tal \cos\tbe+\h_2 \cos\tbe \sin\tal-\h_3 \sin\tbe,\\
\rho_{2}  = -\cos\tga (\h_2 \cos\tal + \h_1 \sin\tal) + (\h_3 \cos\tbe + (-\h_1 \cos\tal + \h_2 \sin\tal) \sin\tbe) \sin\tga ,\\
\rho_{3}  = -\h_3 \cos\tbe \cos\tga + \cos\tga (\h_1 \cos\tal - \h_2 \sin\tal) \sin\tbe - (\h_2 \cos\tal + \h_1 \sin\tal) \sin\tga.
\end{array}
\end{equation*}
The question of integrability of Hamiltonian systems was actively studied by V.I.~Arnold~\cite{3}. Our research continues his study and examines an important example of integrable system.
\end{rem}

\section{Expression for $\displaystyle u_k$, $k\neq 3$ via $u_3$ and the initial values}\label{sec:2}
Let $T>0$, $g \in C(0,T)$. If $g$ is unbounded, assume existence of the integral $\int_0^T g(t) \dd t$. Denote
$$G(t) = \int_0^t g(\tau) \dd \tau.$$
It is known~\cite[ch.~1,~par.~3]{5}, that under such assumptions the Cauchy problem $\dot{y}(t) = g(t) y(t)$, $y(0) = y_0$ has a unique solution 
$y \in C[0,T] \cup \mathcal{D}(0,T)$ given by $y(t) = y_0 \, \eexp(G(t))$. 

Similarly, under the same assumptions the Cauchy problem
\begin{equation}
\label{eq:Cauchysyst}
\begin{cases}
\dot{v}(t) = g(t) w(t), \quad v(0) = v_0, \\
\dot{w}(t) = g(t) v(t), \quad w(0) = w_0
\end{cases}
\end{equation}
has a unique solution $(v,w)$ given by
\begin{equation}
\label{eq:vwsol}
\begin{array}{c}
v(t) = \frac{v_0 + w_0}{2} \eexp\left(G\left(t\right)\right) +  \frac{v_0 - w_0}{2} \eexp\left(-G\left(t\right)\right),\\
w(t) = \frac{v_0 + w_0}{2} \eexp\left(G\left(t\right)\right) -  \frac{v_0 - w_0}{2} \eexp\left(-G\left(t\right)\right).
\end{array}
\end{equation}

Notice that the first and the fifth equations of system~(\ref{eq:extcontsys}) can be written in form~(\ref{eq:Cauchysyst}), where $g(t) = -u_3(t)$, and the second and the fourth equations of system~(\ref{eq:extcontsys}) can be written in form~(\ref{eq:Cauchysyst}), where $g(t) = u_3(t)$. Thus, denoting
\begin{equation}
\label{eq:Ucap}
U(t) = \int_0^t u_3(\tau) \dd \tau
\end{equation} 
and using~(\ref{eq:vwsol}), we express $u_1$, $u_2$, $u_4$, $u_5$ via integral~(\ref{eq:Ucap}) and the initial values
\begin{equation}
\label{eq:extcontrolsuk}
\begin{array}{c}
u_1(t) = \frac{u_1(0) + u_5(0)}{2} \eexp\left(-U\left(t\right)\right) +  \frac{u_1(0) - u_5(0)}{2} \eexp\left(U\left(t\right)\right),\\
u_2(t) = \frac{u_2(0) + u_4(0)}{2} \eexp\left(U\left(t\right)\right) +  \frac{u_2(0) - u_4(0)}{2} \eexp\left(-U\left(t\right)\right),\\
u_4(t) = \frac{u_2(0) + u_4(0)}{2} \eexp\left(U\left(t\right)\right)  -  \frac{u_2(0) - u_4(0)}{2} \eexp\left(-U\left(t\right)\right),\\
u_5(t) = \frac{u_1(0) + u_5(0)}{2} \eexp\left(-U\left(t\right)\right) -  \frac{u_1(0) - u_5(0)}{2} \eexp\left(U\left(t\right)\right).
\end{array}
\end{equation}

\section{Expression for the function $u_3$}\label{sec:3}
It follows from~(\ref{eq:extcontrolsuk}) that
\begin{eqnarray*}
u_1(t) u_5(t) = \left(\frac{u_1(0) + u_5(0)}{2}\right)^2 \eexp\left(-2 U\left(t\right)\right) - \left(\frac{u_1(0) - u_5(0)}{2}\right)^2 \eexp\left(2 U\left(t\right)\right),\nonumber\\
u_2(t) u_4(t) = \left(\frac{u_2(0) + u_4(0)}{2}\right)^2 \eexp\left(2 U\left(t\right)\right) -  \left(\frac{u_2(0) - u_4(0)}{2}\right)^2 \eexp\left(-2 U\left(t\right)\right).\nonumber
\end{eqnarray*}
Therefore,
\begin{equation}
\label{eq:7} 
u_1(t) u_5(t) - u_2(t) u_4(t) = \frac14\left( A \eexp\left(-2 U\left(t\right)\right)  -  B \eexp\left(2 U\left(t\right)\right)\right), 
\end{equation}
where $A = \left(u_1(0) + u_5(0)\right)^2 + \left(u_2(0) - u_4(0)\right)^2$, $B = \left(u_1(0) - u_5(0)\right)^2 + \left(u_2(0) + u_4(0)\right)^2$.

Substitution of~(\ref{eq:7}) in the third equation of system~(\ref{eq:extcontsys}) gives the following second order autonomous differential equation on integral~(\ref{eq:Ucap}):
\begin{equation}
\label{eq:8}
\dot{u}_3(t) = \ddot{U}(t) = \frac{A}{4} \eexp\left(-2 U(t)\right) - \frac{B}{4} \eexp\left(2 U(t)\right).
\end{equation}
There are three possible cases: two special cases ($A=0$ or $B=0$) and the general case $AB \neq 0$ (in this case $A$ and $B$ both are positive). Next we study these 3 cases.

\textbf{I.} $A=0 \Leftrightarrow u_1(0) = -u_5(0), \; u_2(0) = u_4(0)$. Equation~(\ref{eq:8}) becomes
\begin{equation}
\label{eq:9}
\ddot{U}(t) = - B_1 \eexp\left(2 U(t)\right), \quad \text{ where } B_1 = u_4^2(0) + u_5^2(0).
\end{equation}
We aim for a solution that satisfies the initial conditions
\begin{equation}
\label{eq:10}
U(0) = 0, \quad \dot{U}(0) = u_3(0).
\end{equation}
Initial value problem~(\ref{eq:9}), (\ref{eq:10}) can be solved by standard methods. A solution is given by
\begin{equation*}
\begin{array}{l}
U(t) = - \ln\left(\frac12\left[\left(1 + \frac{u_3(0)}{b}\right) \ee^{-b t} + \left(1 - \frac{u_3(0)}{b}\right) \ee^{b t} \right] \right), \; \text{ where } b = \sqrt{u_3^2(0)+u_4^2(0)+u_5^2(0)}.
\end{array}
\end{equation*}
Therefore, we find
\begin{equation}
\label{eq:case1}
u_3(t) = \dot{U}(t) = \frac{\left(b + u_3(0)\right) \ee^{-b t}-\left(b - u_3(0)\right) \ee^{b t}}{\left(1 + \frac{u_3(0)}{b}\right) \ee^{-b t} + \left(1 - \frac{u_3(0)}{b}\right) \ee^{b t}}.
\end{equation}

\textbf{II.} $B=0 \Leftrightarrow u_1(0) = u_5(0), \; u_2(0) = -u_4(0)$. Equation~(\ref{eq:8}) becomes
\begin{equation*}
\ddot{U}(t) = B_1 \eexp\left(-2 U(t)\right).
\end{equation*}
A solution that satisfies initial conditions (\ref{eq:10}) is given by
\begin{equation*}
\begin{array}{l}
U(t) = \ln\left(\frac12\left[\left(1 + \frac{u_3(0)}{b}\right) \ee^{b t} + \left(1 - \frac{u_3(0)}{b}\right) \ee^{-b t} \right] \right).
\end{array}
\end{equation*}
Therefore, we find
\begin{equation}
\label{eq:case2}
u_3(t) = \dot{U}(t) = \frac{\left(b + u_3(0)\right) \ee^{b t}-\left(b - u_3(0)\right) \ee^{-b t}}{\left(1 + \frac{u_3(0)}{b}\right) \ee^{b t} + \left(1 - \frac{u_3(0)}{b}\right) \ee^{-b t}}.
\end{equation}

\textbf{III.} $AB\neq 0 \Rightarrow A>0, \, B>0$. Denote $V = 2 U$, $V_0 = \frac12 \ln \left(\frac{B}{A}\right)$ and rewrite equation~(\ref{eq:8}) as
$$\ddot{V} = \sqrt{AB}\, \frac{\sqrt{A/B} \ee^{-V} - \sqrt{B/A} \ee^{V}}{2} = \sqrt{AB}\, \frac{\ee^{-V-V_0} - \ee^{V+V_0}}{2} = -\sqrt{AB} \sh(V+V_0).$$ 
Next, denoting $y = V+V_0$ we obtain the following Cauchy problem:
\begin{equation}
\label{eq:11}
\ddot{y} = -\sqrt{AB} \sh y, \quad y(0) = \frac12 \ln \left(\frac{B}{A}\right), \; \dot{y}(0) = 2 u_3(0).
\end{equation}
In~\cite{8}, the authors find a solution to problem~(\ref{eq:11}). It leads to
\begin{equation*}
\begin{array}{l}
\displaystyle y(t) = \ln\left(1 + \frac{P^2}{2 \sqrt{A B}}\left(\cn^2\left(\psi_t, k\right) + \frac{1}{k} \cn\left(\psi_t,k\right) \dn\left(\psi_t,k\right) \right)\right),\\
\displaystyle \dot{y}(t) = -P \; \sn\left(\psi_t, k\right),
\end{array}
\end{equation*} 
where
$\displaystyle \psi_t = F\left(p_0, k\right) + \frac{Q}{2} t, \quad
k = \frac{P}{Q}, \qquad p_0 = \begin{cases}
-\arcsin\left(\frac{2 u_3\left(0\right)}{P} \right), \text{ if } B\geq A,\\ 
\pi + \arcsin\left(\frac{2 u_3\left(0\right)}{P} \right), \text{ if } B< A,\end{cases}$\\
with $P = \sqrt{4 u_3^2\left(0\right) + \left(\sqrt{A} - \sqrt{B}\right)^2}$, $\; Q = \sqrt{4 u_3^2\left(0\right) + \left(\sqrt{A} + \sqrt{B}\right)^2}$.

Here, the Jacobi functions $\sn$, $\cn$, $\dn$ and the elliptic integral of the first kind $F$ are used, see~\cite{6}.

Finally, by backward substitutions we express
\begin{equation}
\label{eq:case3}
U\left(t\right)  = \frac{y\left(t\right)}{2} - \frac14 \ln\left(\frac{B}{A}\right), \qquad u_3(t) = \frac{\dot{y}\left(t\right)}{2}.
\end{equation} 
%%%%%%%%%%%%%%%%%%%%%%%%%%%%%%%%%
\section{Conclusion}
Let us summarize results of Sections~\ref{sec:1},~\ref{sec:2},~\ref{sec:3}. The following theorem is proved.
\begin{teo}\label{thm1}
Consider the SR problem in $\SE$. Suppose $u_6(0) = 0$; then vertical part (on extremal controls) of the Hamiltonian system of PMP %system of ODEs on the extremal controls $u_k$, $k \in \{3,4,5\}$  
is given by~(\ref{eq:extcontsys}).\\
The extremal controls $u_4$, $u_5$ are expressed via $U(t) = \int_0^t u_3(\tau) \dd \tau$ and the initial values in~(\ref{eq:extcontrolsuk}).\\
The extremal control $u_3$ is given in terms of the initial values depending on several cases.
For the cases $u_1(0) = \pm u_5(0)$, $u_2(0) = \mp u_4(0)$, we have~(\ref{eq:case1}), ~(\ref{eq:case2}). 
Otherwise, we have~(\ref{eq:case3}).
\end{teo}
In future work, we plan to perform explicit integration of the geodesic equation $\dot{\gamma}(t) = \sum_{i=3}^5 u_i(\tau) \mathcal{A}_i$ as well as study of the general case $\h_6(0) \neq 0$.  
%%%%%%%%%%%%%%%%%%%%%%%%%%%%%%%%%
\section*{Acknowledgments}%in this section you can insert your funding and acknowledgments to somebody in arbitrary form
The authors thank to Prof. Sachkov and Dr. Ardentov for useful discussions. 
%%%%%%%%%%%%%%%%%%%%%%%%%%%%%%%%%

\end{document}